\documentclass[11pt,intlimits, draft,oneside]{amsart}

\usepackage[latin2]{inputenc}
\usepackage{latexsym}
\usepackage{amssymb}
\usepackage{amsthm}

\usepackage{amsmath}

\usepackage{enumerate}

\input{cyracc.def}
\newfont{\cyrfnt}{wncyi10 at 11pt}

\newtheorem{thm}{Theorem}[section]
\newtheorem{prop}[thm]{Proposition}

\newtheorem{cor}[thm]{Corollary}

\theoremstyle{definition}
\newtheorem{defi}[thm]{Definition}
\newtheorem{remark}[thm]{Remark}
\newtheorem{question}[thm]{Question}
\newtheorem{example}[thm]{Example}

\newcommand{\R}{\mathbb{R}}

\newcommand{\C}{\mathbb{C}}
\newcommand{\N}{\mathbb{N}}
\newcommand{\Z}{\mathbb{Z}}

\def\rg{\mathop{\mathrm{rg}}}

\def\LLL{\mathcal{L}}

\begin{document}
\title[Embedding operators into $C_0$-semigroups]
{
Embedding
operators into strongly continuous semigroups}

\author{Tanja Eisner}
\address{Mathematisches Institut,
  Universit\"{a}t T\"{u}bingen \newline
  Auf der Morgenstelle 10, D-72076, T\"{u}bingen, Germany}
\email{talo@fa.uni-tuebingen.de}


\begin{abstract}
We study linear operators $T$ on Banach spaces for which there
exists a $C_0$-semigroup $(T(t))_{t\geq 0}$ such that $T=T(1)$. We
present a necessary condition in terms of the spectral value $0$ and
give classes of examples where this can or cannot be achieved.
\end{abstract}
\maketitle

\begin{center}
\textit{The real understanding involves, I believe, a synthesis of the discrete and continuous ...}

{\small L. Lov\'asz, Discrete and Continuous: two sides of the same?}
\end{center}

\section{Introduction}

%
%

Already the ancient Greek philosophers were asking how time is
flowing:
discrete or conti\-nuous?
\emph{``Into the same river we step and do not step; we are it, and
we are not.''}(Heraclitus of Ephesus, ca. 550--480 B.C.) We refer to
Engel, Nagel \cite[pp. 549--553]{engel/nagel:2000} for some remarks
on this problem.

Even today this is reflected in the mathematical models describing time dependent processes. In the time discrete case one studies a single transformation $\varphi$ and its powers $\varphi^n$, $n\in \N$, on the appropriate state space, while the time continuous case is modelled by a continuous parameter family $\varphi_t$, $t\in\R_+$, of transformations.
Both models are used in the theory of dynamical systems or
stochastic processes and require quite different tools, even if
frequently leading to similar results. So one might look for a
deeper relation between the discrete and the continuous model. This
paper is a contribution to this problem in a particular situation.


We study the question which discrete linear dynamical systems
$\{T^n\}_{n=0}^\infty$ can be embedded into a continuous dynamical
system $(T(t))_{t\geq 0}$. More presicely, we study the following
property.
\begin{defi}
We will say that a linear operator $T$ on a Banach space can be
\textit{embedded} into a $C_0$-semigroup (or, shortly, is
\emph{embeddable}) if there exists a $C_0$-semigroup $T(t)_{t\geq
0}$ on $X$ such that $T=T(1)$.
\end{defi}
\noindent Note that this property implies the existence
of all roots of $T$. 
Note further that $T$ and $(T(t))_{t\geq 0}$ share many properties
such as strong/weak/norm stability, see e.g. Eisner, Farkas, Nagel,
Sereny \cite{EFNS}. However, one never has uniqueness of the
semigroup $T(\cdot)$ since all rescaled semigroups
$(T_n(t)):=(e^{2\pi i n t}T(t))$ satisfy $T_n(1)=T(1)$. We will
concentrate on the existence of a semigroup into which a given
operator can be embedded.

This is a difficult problem which has analogues in other areas of
mathematics like ergodic theory (see e.g. King \cite{king:2000}, de
la Rue, de Sam Lazaro \cite{delarue/desamlazaro:2003}, Stepin,
Eremenko \cite{stepin/eremenko:2004}), stochastics and measure
theory (see e.g. Heyer \cite[Chapter III]{heyer:1977}, Fischer
\cite{fischer:1972}).

It becomes easy if one can construct ``$\ln T$'' by some functional
calculus (see Section 2). In this case the semigroup generated by
$\ln T$ does the job. However, spectral calculus methods cannot be
applied in general since they
require certain conditions on the spectrum of $T$, while, as we see
in Section 2, the spectrum of an embeddable operator can be
arbitrary.

In Section 3 we present a necessary condition for embeddability in
terms of the dimension of $\ker T$ and the codimension of $\rg T$.
Furthermore, we show in Section 4 that for isometries on Hilbert
spaces and their adjoints this condition is also sufficient. In
Section 5 we give more examples of embeddable operators and discuss
the question, how many operators have the embedding property. We
also indicate connections to ergodic theory and the corresponding
results there.



\section{A sufficient condition}

We start with the classical approach using functional calculus. It is based on a sufficient spectral condition
allowing to construct the generator of $T(\cdot)$ as a logarithm of $T$.
\begin{thm}\label{thm:embedding-Dunford}
Assume that $\sigma(T)$ is contained in a simple connected open
domain which does not include $\{0\}$. Then $T$ can be embedded into
a $C_0$-semigroup.
\end{thm}
\begin{proof}
  By the Dunford functional calculus, we can define $A:=\ln T$ as a bounded operator. Then we have
  $T=e^{A}$ and $T$ can be embedded into the semigroup $(T(t))=(e^{tA})$ which is even uniformly continuous.
\end{proof}
\begin{remark}
Note that one can construct $\ln T$ as a bounded operator also in
some other cases. For example, if $X$ is a UMD-space and the Cayley
transform $A:=(1+T)(1-T)^{-1}$ exists and generates a $C_0$-group
with exponential growth less than $1$, then $\ln T$ exists and is
bounded, so $T$ can be embedded into a uniformly continuous
$C_0$-semigroup. For details see e.g. Haase \cite{haase:2007}.
\end{remark}

There are several extensions of Theorem \ref{thm:embedding-Dunford}.
One of them is the following result which allows to construct
semigroups with unbounded generators, hence corresponding to
unbounded logarithms. Recall that an operator is called
\emph{sectorial} if there exists a sector
$\Sigma_\delta=\{z:|\text{arg} z|\leq \delta\}\cup\{0\}$,
$0<\delta<\pi$, such that $\sigma(T)\subset \Sigma_\delta$ and for
every $0<\delta<\omega$ one has $\|R(\lambda,T)\|\leq
\frac{M}{|\lambda|}$ for some $M$ and every $\lambda\notin
\Sigma_\omega$.

\begin{thm}\label{thm:embedding-haase}(See Haase \cite[Prop. 3.1.1 and 3.1.15]{haase:book})
Let $T$ be a bounded sectorial operator with dense range. Then $T$ can be embedded into an analytic $C_0$-semigroup.
\end{thm}
\noindent Since $T$ is embeddable if and only if $cT$ is embeddable for any $0\neq c \in \C$, we see that one can consider
operators with spectrum in a rotated sector $\{z:|\text{arg} z-\varphi|<\delta\}$ as well.


\smallskip

We now show that the spectrum of an embeddable operator can be
arbitrary, hence conditions on the location of the spectrum are not necessary and
 the spectral calculus method works only in particular cases.

\begin{example}\label{ex:embedding-spectrum}
Let $K$ be a compact set in $\C$. Assume first that $0$ is not an
isolated point in $K$ and consider a dense subset
$\{\lambda_n\}_{n=-\infty}^\infty$ of $K\setminus \{0\}$. Consider
further $X:=l^2$ and the multiplication operator defined by
$T(x_1,x_2,\ldots):=(\lambda_1 x_1, \lambda_2 x_2, \ldots)$. Observe
that $\sigma(T)=K$ while $T$ can be embedded into the semigroup
\[
T(t)(x_1,x_2,\ldots):=(e^{t\ln\lambda_1} x_1, e^{t\ln\lambda_2} x_2,
\ldots), \quad t\geq 0,
\]
where we take an arbitrary value of the logarithm for every
$\lambda_j$. Since $T(t)e_j$ is continuous in $t$ for every basis
vector $e_j$ and $(T(t))_{t\geq 0}$ is uniformly bounded on compact
time intervals, it is a $C_0$-semigroup.

Assume now that $\{0\}$ is an isolated point of $K$, hence
$K=\{0\}\cup K_1$ for some compact set $0\notin K_1$. Define
the embeddable operator $T_1$ with $\sigma(T_1)=K_1$ as above. Define
further $T_2=0$ on $X=l^2$. The operator $T_2$ can be embedded into a
nilpotent semigroup on $l^2$. To construct such a semigroup we first
consider $L^2[0,1]$ and the nilpotent shift semigroup on it. Since
all separable infinite-dimensional Hilbert spaces are isomorphic, it
corresponds to a nilpotent $C_0$-semigroup on $l^2$ in which we
embedd $T_2$. The operator $T_1\oplus T_2$ on $X^2$ is an example of an embeddable operator with spectrum $K$.

\end{example}
\begin{remark}
Using an analogous method, one can construct for every compact set $K$
a non-embeddable operator $T$ with $\sigma(T)=K\cup\{0\}$. We just take the direct sum
of the operator constructed in the first part of Example
\ref{ex:embedding-spectrum} with the zero operator on the one-dimensional
space. This sum is not embeddable by Theorem
\ref{thm:embedding-necessary-cond} below.
\end{remark}

\section{A necessary condition}

In this section we present a simple necessary condition for being
embeddable using information on the spectral value $0$.

\begin{thm}\label{thm:embedding-necessary-cond}
Let $X$ be a Banach space and $T\in\LLL(X)$.
If $T$ can be embedded into a $C_0$-semigroup, then $\dim(\ker T)$ and $\emph{codim}(\rg T)$ are equal to zero or to infinity.
\end{thm}
\noindent In other words, operators with $0< \dim(\ker T)<\infty$
or $0< \text{codim}(\rg T)<\infty$ cannot be embedded into a $C_0$-semigroup.
\begin{proof}
Assume that $0<\dim(\ker
T)<\infty$ holds and $T=T(1)$ for some $C_0$-semigroup $T(\cdot)$ on $X$. Since
$T$ is not injective, so is its square root
$T\left(\frac{1}{2}\right)$. Analogously, every $T\left(\frac{1}{2^n}\right)$ is not
injective. Take $x_n\in \ker T\left(\frac{1}{2^n}\right)$ with $\|x_n\|=1$ for
every $n\in\N$. Then we have $\{x_n\}_{n=1}^\infty\subset \ker T$
which is finite-dimensional. Therefore there ixists a subsequence $\{x_{n_k}\}_{k=1}^\infty$
converging to some $x_0$ with $\|x_0\|=1$.

Since $\ker T\left(\frac{1}{2^{n+1}}\right)\subset \ker T\left(\frac{1}{2^n}\right)$,
$n\in\N$, we have $T\left(\frac{1}{2^n}\right)x_0=0$ for every $n\in\N$ which
contradicts the strong continuity of $T(\cdot)$.

Assume now that $0< \text{codim}(\rg T)$$=\dim(\ker
T')<\infty$ holds. If $T$ is embedded into a $C_0$-semigroup $T(\cdot)$, then $T'$ is embedded into the adjoint semigroup $T'(\cdot)$ on $X'$ which is weak* continuous. The same arguments as above lead to a contradiction to the weak* continuity.
\end{proof}
A direct corollary is the following.
\begin{cor}
Non-bijective Fredholm-operators are not embeddable.
\end{cor}
\begin{remark}\label{remark:embedding-fore-possib}
As we will see in Section \ref{section-embedding-isometry}, all fore
possibilities given in Theorem \ref{thm:embedding-necessary-cond}
(dimension $0$ or $\infty$ of the kernel and codimension $0$ or
$\infty$ of the range) can appear for operators with the embedding property.
The examples will be unitary operators, the left and right shifts on
$l^2(Y)$ for an infinite-dimensional Hilbert space $Y$ and their
direct sum, respectively.
\end{remark}

In the following we ask the converse question.
\begin{question}\label{question:embedding-back}
For which (classes of) operators the necessary condition given in
Theorem \ref{thm:embedding-necessary-cond} is also sufficient for embeddability?
\end{question}
%


\section{Embedding isometries on Hilbert spaces}\label{section-embedding-isometry}

In this section we characterise the embedding property for
isometries on Hilbert spaces. Note that the spectrum of a
non-invertible isometry is the unit disc, and hence
the spectral calculus method is not applicable.

We first recall that unitary operators are embeddable.

\begin{prop}\label{prop:embedding-unitary}
Every unitary operator $U$ on a Hilbert space $H$ can be embedded into a unitary $C_0$-group with bounded generator.
\end{prop}
\begin{proof}
By the spectral theorem, see e.g. Halmos \cite{halmos:1963}, we can
assume that $U$ has a cyclic vector, $H=L^2(\Gamma, \mu)$ for the unit circle $\Gamma$ and some Borel measure $\mu$ and
\begin{equation*}
(Uf)(z)=zf(z), \quad z\in\Gamma, \ f\in L^2(\Gamma, \mu).
\end{equation*}
Define now
\begin{equation*}
(\tilde{U}(t)f)(e^{i\varphi}):=e^{it\varphi}f(e^{i\varphi}), \quad f\in L^2(\Gamma, \mu),
\end{equation*}
for $\varphi\in [0,2\pi)$ and $t\in\R$. The operators $\tilde{U}(t)$
form a unitary $C_0$-group which is even uniformly continuous and satisfies $U(1)=U$.
\end{proof}

To treat the general case we need the following structure
theorem for isometries on Hilbert spaces.
\begin{thm}\label{thm:Wold}({\it Wold decomposition}, see Sz.-Nagy, Foia{\c{s}} \cite[Theorem 1.1]{sznagy/foias}.)
Let $V$ be an isometry on a Hilbert space $H$. Then $H$ can be decomposed into an orthogonal sum $H=H_0 \oplus H_1$ of $V$-invariant subspaces such that the restriction of $V$ on $H_0$ is unitary and the restriction of $V$ on $H_1$ is a unilateral shift. More precisely, for $Y:=(\rg V)^\perp \subset H_1$ one has $V^n Y\perp V^m Y$ for all $n\neq m\in\N$ and $H_1=\oplus_{n=1}^\infty V^n Y$.
\end{thm}

So, by Wold's decomposition and Proposition \ref{prop:embedding-unitary}, the question of embedding an isometry restricts to embedding a right shift on the space $l^2(Y)$ for some Hilbert space $Y$.

The following result shows when this can be achieved.
%
\begin{prop}\label{prop:embedding-shift-infinite}
Let $S$ be the right shift on $l^2(Y)$ for an infinite-dimensional Hilbert space $Y$. Then
$S$ can be embedded into an isometric $C_0$-semigroup.
\end{prop}
\begin{proof}
Let $T$ be the right shift on $l^2(Y)$, i.e.,
\begin{equation*}
T(x_1,x_2,x_3,\ldots):=(0,x_1,x_2,\ldots).
  \end{equation*}
By general Hilbert space theory,
$Y$ is (unitarily) isomorphic to the space $L^2([0,1],Y)$, hence
there is a unitary operator $J:l^2(Y)\to l^2(L^2([0,1],
Y))$ such that $JTJ^{-1}$ is again the right shift operator on
$l^2(L^2([0,1],Y))$.
We now observe that $l^2(L^2([0,1],Y))$ can be identified with $L^2(\R_+,Y)$ by
\begin{equation*}
(f_1,f_2,f_3,\ldots)\mapsto (s\mapsto f_n(s-n),\ s\in [n,n+1]).
  \end{equation*}
(Note that the above identification is unitary.)
Under this identification the right shift on $l^2(L^2([0,1],Y))$
corresponds to the operator
\begin{equation*}
(Sf)(s):=
\begin{cases}
f(s-1),\quad &s\geq 1,\\
0, \quad &s\in[0,1)
  \end{cases}
  \end{equation*}
on $L^2(\R_+,Y)$ which clearly can be embedded into the right shift semigroup on
$L^2(\R_+,Y)$. Going back we see that our original operator $T$ can be
embedded into an isometric $C_0$-semigroup.
%
\end{proof}

%
%
%

We are now ready to give a complete answer to Question
\ref{question:embedding-back} for isometries on Hilbert spaces.
\begin{thm}
An isometry $V$ on a Hilbert space can be embedded into a $C_0$-semigroup if and only if $V$ is unitary or $\emph{codim}(\rg V)=\infty$.
\end{thm}
\begin{proof}
Let $V$ be isometric on a Hilbert space $H$. By the Wold decomposition we have the orthogonal decomposition $H=H_1\oplus H_2$ in two invariant subspaces such that $V|_{H_0}$ is unitary and $V|_{H_1}$ is unitarily equivalent to the right shift on $l^2(Y)$ for $Y:=(\rg V)^\perp$.
By Proposition \ref{prop:embedding-unitary}, we can embedd $V|_{H_0}$ into a unitary $C_0$-group.
On the other hand, by Proposition \ref{prop:embedding-shift-infinite} and Theorem \ref{thm:embedding-necessary-cond} we can embed $V|_{H_1}$ if and only if $Y=0$ or $\dim Y=\infty$.
\end{proof}

Since an operator on a Hilbert space can be embedded into a
$C_0$-semigroup if and only if its adjoint can, we directly obtain the following.
\begin{cor}
Let $T$ be an operator on a Hilbert space with isometric adjoint.
Then $T$ can be embedded into a $C_0$-semigroup if and only if $T$
is injective or $\dim (\ker T)=\infty$.
\end{cor}
\noindent For example, the left shift on $l^2(Y)$ for a Hilbert space $Y$ has the embedding propery if and only if $\dim(Y)=\infty$.

\smallskip

\section{More examples: normal and compact operators, a category result}

In this section we give more examples of operators possessing the embedding property.
\begin{prop}
Let $T$ be a normal operator on a Hilbert space. Then $T$ can be embedded into a normal $C_0$-semigroup.
\end{prop}
\begin{proof}
By the spectral theorem, we can assume that $T$ has a cyclic vector
and has the form $(Tf)(z)=zf(z)$ on $L^2(\Omega,\mu)$ for some
$\Omega\subset \{z:|z|\leq 1\}$ and some measure $\mu$. Define now
\[
(T(t)f)(re^{i\varphi}):= r^{t} e^{it\varphi}f(re^{i\varphi}), \quad
f\in L^2(\Omega,\mu), \ 0<r<1, \ \varphi\in[0,2\pi),
\]
\noindent and $(T(t)f)(0):=0$. Then $T(\cdot)$ is a contractive
normal $C_0$-semigroup satisfying $T(1)=T$.
\end{proof}
\begin{remark}
Clearly every operator similar to a normal operator is also embeddable. For example, Sz.-Nagy \cite{sznagy:1947} showed that a bijective operator $T$ on a Hilbert space satisfying $\sup_{j\in\Z} \|T^j\|<\infty$ is similar to a unitary operator. Hence such operators can be embedded into a $C_0$-group by Proposition \ref{prop:embedding-unitary}. We refer also to van Casteren \cite{vancasteren:1983} for conditions to be similar to a self-adjoint operator.
We finally refer to Benamara, Nikolski \cite{benamara/nikolski:1999} for a characterisation of operators similar to a normal operator as well as for the history of the similarity problem.
\end{remark}

To give more examples of embeddable operators we consider compact
operators on Banach spaces.
\begin{prop}
Let $T$ be an injective compact operator with dense range on a Banach space $X$. Then
$T$ can be embedded into a $C_0$-semigroup.
\end{prop}
\begin{proof}
By the spectral theorem for compact operators, we can decompose
$X=X_0\oplus X_1 \oplus X_2 \oplus \ldots$ and correspondingly
$T=T_0\oplus T_1 \oplus T_2\oplus \ldots$ with finite-dimensional
invariant subspaces $X_1,X_2,\ldots$ such that $\sigma(T_0)=\{0\}$
and $\sigma(T_j)=\{\lambda_j\}$ for $\lambda_j\neq 0$. Every $T_j$
can be embedded by spectral calculus. Since $r(T_0)=0$, the operator $T_0$ is sectorial by the Neumann series representation for the resolvent and hence can be embedded as well by Theorem \ref{thm:embedding-haase}.
So the original operator $T$ can be embedded into a $C_0$-semigroup.
\end{proof}
\noindent As an example we see that the Volterra operator
$(Vf)(\tau)=\int_0^\tau f(s)ds$ on $X=C[0,1]$ has the embedding
property.
\begin{remark}
It is not clear whether a compact operator with dense (or
infinite-codimensional)
range and infinite-dimensional kernel can always be embedded into a $C_0$-semigroup.
(Note that finite-dimensional kernel and cofinite-dimensional range is
excluded by Theorem \ref{thm:embedding-necessary-cond}.) Furthermore,
it is desirable to construct a non-trivial operator (in the
sense of Theorem \ref{thm:embedding-necessary-cond}) which cannot be embedded into a $C_0$-semigroup.
\end{remark}

\vspace{0.05cm}

As a final remark we mention that for separable Hilbert spaces a ``typical'' contraction is embeddable.
\begin{thm}(see Eisner \cite{eisner:2008})
The set of all embeddable contractions on a separable
infinite-dimensional Hilbert space form a residual set for the weak
operator topology.
\end{thm}
\noindent To prove this fact, one first shows that unitary operators form a
residual set (see Eisner \cite[Theorem 3.3]{eisner:2008}). Then Proposition \ref{prop:embedding-unitary}
finishes the argument.

\begin{remark}
This result is an operator-theoretical counterpart to a recent result
of de la Rue and de Sam Lasaro \cite{delarue/desamlazaro:2003} from ergodic theory stating that a
``typical'' measure preserving transformation can be embedded into a
measure preserving flow. We refer to Stepin, Eremenko \cite{stepin/eremenko:2004} for further results and references.

In particular, a ``typical'' contraction on a separable
infinite-dimensional Hilbert space has roots of all order, which is an operator-theoretical analogue of a result of King \cite{king:2000} in ergodic theory.
\end{remark}

\smallskip

\noindent {\bf Acknowledgement.} The author is very grateful to
Rainer Nagel for interesting discussions.

\parindent0pt

\end{document}